\documentclass[final,3p,times]{elsarticle}

\usepackage{amssymb}
\usepackage{graphicx,setspace}
\usepackage{psfrag}
\usepackage{amsfonts}
\expandafter\let\csname equation*\endcsname\relax
\expandafter\let\csname endequation*\endcsname\relax
\usepackage{amsmath}
\usepackage{amssymb, amsthm}
\usepackage{mathrsfs}
\usepackage{epstopdf}
\usepackage{float}
\usepackage{lineno,hyperref}
\usepackage{color}
\usepackage{subfigure}
\usepackage{geometry}
\usepackage{comment}
\usepackage{booktabs}
\usepackage{caption}
\usepackage{hyperref}
\usepackage{amsmath}
\usepackage{dsfont} % it is the package of using the "\mathds{} "
\usepackage{amssymb} % it is the package of using the "\mathbb{}"
\usepackage{graphicx,color}
\usepackage{extarrows}

\usepackage{graphicx}

\usepackage{epstopdf}
\journal{arXiv}

\begin{document}
\newtheorem{definition}{Definition}[section]
\newtheorem{lemma}{Lemma}[section]
\newtheorem{remark}{Remark}[section]
\newtheorem{theorem}{Theorem}[section]
\newtheorem{proposition}{Proposition}
\newtheorem{assumption}{Assumption}
\newtheorem{example}{Example}
\newtheorem{corollary}{Corollary}[section]
\def\ep{\varepsilon}
\def\Rn{\mathbb{R}^{n}}
\def\Rm{\mathbb{R}^{m}}
\def\E{\mathbb{E}}
\def\hte{\hat\theta}
%\numberwithin{theorem}{section}
%\numberwithin{definition}{section}
\renewcommand{\theequation}{\thesection.\arabic{equation}}
\begin{frontmatter}



\title{Solow system driven by $\alpha$-stable L\'evy process}

\author{Yiren Wang\fnref{addr1,addr2}}\ead{252221004@stu.gbu.edu.cn}
\author{Shenglan Yuan\corref{cor1}\fnref{addr1}}
\ead{shenglanyuan@gbu.edu.cn}\cortext[cor1]{Corresponding author}

\address[addr1]{\rm Department of Mathematics, School of Sciences, Great Bay University, Dongguan 523000, China }
\address[addr2]{\rm Department of Mathematics, College of Science, Southern University of Science and Technology,  Shenzhen 518055,
 China}

\begin{abstract}
	This paper develops and empirically implements a Solow-type growth
		model driven by symmetric $\alpha$-stable L\'evy shocks with
		time-varying capital elasticity. We extend the framework of Brannan
		(2019) by replacing the Gaussian Ornstein-Uhlenbeck process with an
		$\alpha$-stable L\'evy process ($\alpha\in(1,2)$), thereby capturing
		three stylized facts of severe macroeconomic fluctuations: heavy-tailed
		distributions, jump discontinuities, and infinite variance.
		Theoretically, we derive the stationary distribution of the capital
		deviation process, obtain its conditional characteristic function in
		closed form, and provide an integral representation that explicitly
		reveals a dual mean-reversion structure separating investment gestation
		lags from endogenous feedback. Methodologically, we design an estimation
		strategy based solely on well-defined objective functions---mean absolute
		prediction error and mean matching---that respects the probabilistic
		properties of L\'evy-driven data and circumvents the non-existence of
		variance. Empirically, we apply the framework to Argentine quarterly
		data from 2004 to 2023, with time-varying capital elasticity calibrated
		from Penn World Table labor shares. Our estimates show that the L\'evy
		specification delivers structural parameters substantially closer to
		external PWT benchmarks than the Gaussian Ornstein-Uhlenbeck counterpart and
		substantially improves crisis-period tracking without sacrificing
		performance in tranquil periods. Cross-country evidence from Colombia and the United States confirms that the quarterly capital adjustment speed \(\eta \approx 0.05\)---a genuine structural parameter reflecting the physical speed of fixed capital formation---exhibits striking stability across vastly different volatility regimes, whereas the Gaussian Ornstein-Uhlenbeck benchmark systematically inflates this parameter by a factor of 7 to 10.
		Robustness checks across tail index specifications demonstrate that
		the L\'evy framework consistently outperforms the Ornstein-Uhlenbeck benchmark for a
		broad range of empirically relevant tail indices. These findings
		establish the L\'evy specification as a robust generalization of the
		Gaussian benchmark, offering a more credible tool for forecasting and
		structural parameter estimation in both emerging and advanced economies.
\end{abstract}


\end{frontmatter}


\section{Introduction}
	Economic growth and business cycles represent core themes of macroeconomic research.
	The Solow growth model \cite{solow1956contribution} provides the fundamental framework for understanding long-run economic growth,
	describing the evolution of capital stock per effective unit of labor \(k(t)\) through the dynamic equation
	\[
	\frac{dk}{dt} = s f(k) - (n + g + \delta)k,
	\]
	where \(s\) is a constant saving rate, \(f(k)\) stands for the production function, \(n\) represents the labor growth rate,
	\(g\) denotes the rate of technological progress, and \(\delta\) is the capital depreciation rate.
	Under standard assumptions, the economy converges to a balanced growth path (BGP) equilibrium \(k^*\)
	where investment exactly offsets effective depreciation  $(n + g + \delta)k$.
	Although the Solow model successfully explains long-run convergence and steady-state growth,
	it is essentially deterministic and cannot fully explain endogenous cyclical fluctuations
	or realistic discontinuous economic shocks.
	
	Brannan \cite{brannan2019natural} significantly extends the Solow framework by introducing three key modifications:
	endogenous saving behavior, time delays in investment realization, and stochastic disturbances.
	First, the constant saving rate is replaced by a sigmoidal function of capital deviations from the BGP,
	which introduces procyclical positive feedback and generates endogenous persistence.
	Second, investment gestation lags are modeled using an exponentially weighted moving average,
	leading to a fast-slow dynamical system with capital as the slow variable.
	Third, correlated shocks are introduced via an Ornstein-Uhlenbeck (OU) process driven by a Brownian motion.

	However, the Brownian motion framework only describes continuous Gaussian shocks,
	while real economic shocks, especially in emerging and small open economies, are often
	discontinuous and heavy-tailed.
	Sudden crises, sharp shifts in global demand, or extreme policy shocks
	cannot be accurately captured by diffusion processes. The limitations of Gaussian shocks in capturing macroeconomic tail
	risks have motivated a growing literature on fat-tailed
	macroeconomic models, including fat-tailed DSGE specifications
	\cite{dave2025fattailed}. Our $\alpha$-stable L\'evy framework
	belongs to the same broad ``fat in, fat out'' category---where
	non-Gaussian exogenous disturbances directly generate heavy-tailed
	outcomes---but operates within a Solow-type dynamical system rather
	than a micro-founded general equilibrium setting.
	
	To address this limitation, this work extends the Solow-type models
	by replacing Brownian motion with an \(\alpha\)-stable L\'evy process where $1 < \alpha < 2$. The $\alpha$-stable L\'evy process naturally incorporates three stylized facts of severe macroeconomic fluctuations: (i) discontinuities (sudden jumps), (ii) heavy tails (extreme events occur with non-negligible probability), and (iii) infinite variance (reflecting profound uncertainty during crises).
	The use of \(\alpha\)-stable L\'evy processes to drive slow-fast dynamical
	systems has been rigorously developed in the stochastic bifurcation
	literature \cite{yuan2021stochastic}, establishing stochastic equilibrium
	states and bifurcation phenomena for two-time-scale systems driven by
	\(\alpha\)-stable L\'evy noise, providing the theoretical foundation for
	our dynamical system (2).

	Our theoretical contribution is threefold. First, we derive the stationary solution of the capital deviation process and rigorously establish its strict stationarity using characteristic functions. Second, we obtain the conditional characteristic function of future capital deviations given the current state, which fully characterizes the conditional distribution and admits closed-form expressions for conditional moments. Third, we provide a tractable integral representation of the capital deviation as a linear combination of two independent OU processes, explicitly revealing the dual mean-reversion structure driven by investment gestation lags and endogenous feedback. These results provide the probabilistic foundation for all subsequent empirical analysis.
	
	Our empirical contribution is equally important. We apply the proposed framework to Argentina, a typical emerging economy that has experienced recurrent extreme events over the past two decades. Recognizing that production structures in volatile economies are unlikely to remain constant, we allow the capital elasticity of output to vary over time, calibrated through external labor-share data, to accommodate structural breaks. Furthermore, because the $\alpha$-stable shocks lack finite variance, conventional moment-based estimation is infeasible. We therefore design a coherent estimation strategy that relies exclusively on well-defined quantities, respecting the probabilistic properties of the driving process. When applied to Argentine quarterly data, our model yields two main findings. First, the estimated key structural parameters are substantially closer to external benchmarks than those obtained from the standard Gaussian OU specification, confirming the empirical relevance of our modifications. Second, the model's predictions track the observed capital dynamics reasonably well across different sub-periods, with accuracy varying systematically with macroeconomic volatility---a pattern consistent with the inherent non-Gaussian nature of the shocks.
	
	The remainder of this paper is organized as follows. Section 2 develops the continuous-time Solow model driven by $\alpha$-stable L\'evy shocks, derives the stationary solution and integral representation of the capital deviation process, and presents the characteristic functions for the state variables and their future increments. Section 3 establishes the conditional characteristic function and conditional expectation, and provides the trajectory simulation scheme. Section 4 constructs a discrete-time approximation for empirical implementation, details the estimation strategy, presents the empirical results for Argentina, and conducts model comparisons. Section 5 presents robustness checks and cross-country extensions. Section 6 concludes with a discussion of main contributions, policy implications, and avenues for future research.
	
\section{Model setup}\label{Model Setup}
	In this section, we develop our full dynamic model under $\alpha$-stable L\'evy shocks, building on the scaled Solow framework extended by \cite{brannan2019natural}. Our key modification replaces the standard Brownian motion driving the stochastic shocks with an $\alpha$-stable L\'evy process. We derive the resulting stochastic differential system governing the dynamics of capital deviations $u(\tau)$, which forms the basis for our subsequent analysis.
	
	\subsection{The scaled slow-time dynamical system}
	Following \cite{brannan2019natural}, we work in the dimensionless slow time scale $\tau = \tilde{\eta} t$, where $\tilde{\eta} = n + g + \delta$. Let $k$ denote capital per effective worker, scaled by its BGP value, and let $z$ be the exponentially weighted capital stock that captures investment gestation lags.
	
Brannan	\cite{brannan2019natural} incorporates stochastic shocks modeled as a Brownian-motion-driven OU process and linearizes the model around the BGP to obtain the baseline dynamical system:
\begin{equation}\label{eq:baseline}
		\begin{cases}
			\dfrac{du}{d\tau} &= -r(\gamma) u - (\beta-1) Q(\tau), \\[2pt]
			dQ(\tau) &= -a Q(\tau) d\tau + \sigma dW(\tau),
		\end{cases}
\end{equation}
where $a>0$ represents the mean reversion parameter and $\sigma>0$ is the scale parameter.
	
We preserve the entire structural form of this baseline model and only generalize the shock specification. Specifically, we replace the standard Brownian motion $W(\tau)$ with a strictly $\alpha$-stable L\'evy process $L_\alpha(\tau)$ with $\alpha\in(1,2)$. This yields the following dynamical system:
\begin{equation}\label{eq:extended}
\begin{cases}
			\dfrac{du}{d\tau} &= -r(\gamma) u - (\beta-1) Q(\tau), \\[4pt]
			dQ(\tau) &= -a Q(\tau) d\tau + \sigma dL_\alpha(\tau).
\end{cases}
\end{equation}
	
This introduces a novel dynamic characterized by jump discontinuities, heavy tails, and infinite variance. By explicitly capturing these non-Gaussian features, our framework addresses a critical void in the existing literature, enabling robust analysis of extreme shocks and tail risks that lie beyond the scope of traditional linear models. Thus the system \eqref{eq:extended} serves as the foundation for all subsequent analysis.

\subsection{Stationary solutions and integral representations}
The noise process \(Q(\tau)\) itself admits the OU representation \cite{applebaum2009levy, protter2012stochastic}:
	\begin{equation}
		Q(\zeta) = \sigma \int_{-\infty}^{\zeta} e^{-a(\zeta-s)} dL_{\alpha}(s).
		\label{eq:Q_inter_result}
	\end{equation}
	The solution to the capital deviation equation is
	\[
	u(\tau) = (1-\beta) \int_{-\infty}^{\tau} e^{-r(\gamma)(\tau-\zeta)} Q(\zeta) d\zeta.
	\]
	Substituting the steady-state solution for \(Q(\zeta)\) from \eqref{eq:Q_inter_result} into the expression for \(u(\tau)\) yields a double integral:
	\begin{equation}
		u(\tau) = (1-\beta)\sigma \int_{-\infty}^{\tau} e^{-r(\gamma)(\tau-\zeta)} \left( \int_{-\infty}^{\zeta} e^{-a(\zeta-s)} dL_{\alpha}(s) \right) d\zeta.
		\label{eq:u_inter_result}
	\end{equation}
	Now we verify the stationarity of \(Q(\tau)\) and \(u(\tau)\).
\subsubsection{Stationarity of \(Q(\tau)\)}
For any \(h \in \mathbb{R}\), we compute \(Q(\tau+h)\):
	\[
	\begin{aligned}
		Q(\tau+h) &= \sigma \int_{-\infty}^{\tau+h} e^{-a(\tau+h-s)} dL_{\alpha}(s) \\
		&= \sigma \int_{-\infty}^{\tau} e^{-a(\tau-s')} dL_{\alpha}(s'+h) \quad (\text{by the change of variable } s' = s-h) \\
		&\stackrel{d}{=} \sigma \int_{-\infty}^{\tau} e^{-a(\tau-s')} dL_{\alpha}(s') = Q(\tau),
	\end{aligned}
	\]
	where the last equality follows from the stationary increment property of \(L_{\alpha}(t)\).
	The proof uses the standard change-of-variable technique for stochastic integrals, as applied in the context of slow-time systems by Yuan et al. \cite{yuan2019slow}.
	
\subsubsection{Stationarity of \(u(\tau)\)}
	
For any \(h \in \mathbb{R}\), we compute \(u(\tau+h)\):
	\[
	\begin{aligned}
		u(\tau+h) &= (1-\beta)\sigma \int_{-\infty}^{\tau+h} e^{-r(\gamma)(\tau+h-\zeta)} Q(\zeta) d\zeta \\
		&= (1-\beta)\sigma \int_{-\infty}^{\tau} e^{-r(\gamma)(\tau-\zeta')} Q(\zeta'+h) d\zeta' \quad (\text{by the change of variable } \zeta' = \zeta-h) \\
		&\stackrel{d}{=} (1-\beta)\sigma \int_{-\infty}^{\tau} e^{-r(\gamma)(\tau-\zeta')} Q(\zeta') d\zeta' = u(\tau),
	\end{aligned}
	\]
	where the last equality follows from the strict stationarity of \(Q(\tau)\).
	The proof structure mirrors the change-of-variable approach for verifying stationarity in L\'evy-driven systems, as introduced in \cite{yuan2019slow}.

\subsubsection{ Steady-state approximation}
The $\alpha$-stable L\'evy process $L_\alpha(\tau)$ with
	$\alpha\in(1,2)$ exhibits power-law tails and jump discontinuities,
	and its infinitesimal increments are not necessarily absolutely
	integrable. \cite{yuan2022stochastic} demonstrate that the
	stability index $\alpha$ governs the jump intensity and tail
	decay of stochastic trajectories, with lower values producing
	more frequent extreme events. In order to justify the
	interchange of integration order under Fubini's theorem, we
	first establish the stationarity of $Q(\tau)$ and $u(\tau)$.
	Since $L_{\alpha}(\tau)$ is a process with stationary increments, the OU process $Q(\tau)$ is strictly stationary.
	Moreover, $u(\tau)$ given in equation \eqref{eq:u_result} is an exponentially weighted integral of the stationary process $Q(\tau)$, so $u(\tau)$ is also strictly stationary.
	Under stationarity, contributions from the distant past are strongly dampened by the exponential kernels $e^{-a(\zeta-s)}$ and $e^{-r(\gamma)(\tau-\zeta)}$.
	Thus the integrand is absolutely integrable, and the stationary approximation is valid.

	We now interchange the order of integration over the domain $-\infty < s < \zeta < \tau$. Starting from equation \eqref{eq:u_inter_result},
	we obtain
	\begin{align}
		u(\tau) &= (1-\beta)\sigma \int_{-\infty}^{\tau} \left( \int_{s}^{\tau} e^{-r(\gamma)(\tau-\zeta)}e^{-a(\zeta-s)}d\zeta \right) dL_{\alpha}(s), \nonumber \\
		&= (1-\beta)\sigma \int_{-\infty}^{\tau} e^{-a(\tau-s)} \left( \int_{s}^{\tau} e^{-(r(\gamma)-a)(\tau-\zeta)}d\zeta \right) dL_{\alpha}(s), \nonumber \\
		&= (1-\beta)\sigma \int_{-\infty}^{\tau} e^{-a(\tau-s)} \frac{1-e^{-(r(\gamma)-a)(\tau-s)}}{r(\gamma)-a} dL_{\alpha}(s), \nonumber \\
		&= \frac{(1-\beta)\sigma}{a-r(\gamma)} \left( \int_{-\infty}^{\tau} e^{-a(\tau-s)}dL_{\alpha}(s) - \int_{-\infty}^{\tau} e^{-r(\gamma)(\tau-s)}dL_{\alpha}(s) \right).
		\label{eq:u_result}
	\end{align}
This shows that $u(\tau)$ is a linear combination of two $\alpha$-stable OU processes with decay rates $a$ and $r(\gamma)$, and remains strictly stationary.
	
\subsection{Characteristic functions}
	
\subsubsection{Characteristic functions of $Q(\tau)$}\label{subsec:char_func_qu}
	
	In this section, we derive the characteristic functions associated with the stationary processes $Q(\tau)$ and $u(\tau)$.
	Given that both processes are driven by a symmetric $\alpha$-stable L\'{e}vy process $L_{\alpha}(\tau)$,
	their probabilistic properties are fully characterized by characteristic functions.

We first consider the stationary OU process $Q(\tau)$ defined as
	\begin{equation}\label{eq:Q_OU_def}
		Q(\tau) = \sigma \int_{-\infty}^{\tau} e^{-a(\tau-s)} dL_{\alpha}(s).
	\end{equation}
	By the L\'{e}vy-Khintchine representation and the stationary increment property of L\'{e}vy processes,
	for a symmetric $\alpha$-stable L\'evy process $L_{\alpha}(s)$, the characteristic function \cite{protter2012stochastic} of a stochastic integral $X = \int_{s_1}^{s_2} f(s) dL_{\alpha}(s)$ is given by
	\[
	\mathbb{E}\left[ e^{i\theta X} \right] = \exp\left( - C_{\alpha} |\theta|^{\alpha} \int_{s_1}^{s_2} |f(s)|^{\alpha} ds \right),
	\]
	where $C_{\alpha}$ is the scale constant of the stable symmetric distribution $\alpha$-stable distribution.

Thus the characteristic function of $Q(\tau)$ takes the form
	\begin{equation}\label{eq:char_Q_1}
		\E\left[ e^{i\theta Q(\tau)} \right]
		= \exp\left( -C_{\alpha}\sigma^{\alpha} |\theta|^{\alpha} \int_{0}^{\infty} e^{-a\alpha t} dt \right),
	\end{equation}
	where $C_{\alpha}>0$ denotes the scale parameter of the symmetric $\alpha$-stable distribution and $\theta\in\mathbb{R}$.
	Evaluating the deterministic integral yields
	\[
	\int_{0}^{\infty} e^{-a\alpha t} dt = \frac{1}{a\alpha}.
	\]
	Substituting this result into \eqref{eq:char_Q_1}, we obtain
	\begin{equation}\label{eq:char_Q_final}
		\E\left[ e^{i\theta Q(\tau)} \right]
		= \exp\left( -\frac{C_{\alpha}\sigma^{\alpha}}{a\alpha} |\theta|^{\alpha} \right),
	\end{equation}
	which confirms that $Q(\tau)$ follows a strictly stationary symmetric $\alpha$-stable distribution.

\subsubsection{Characteristic functions of $u(\tau)$}
	Next, we derive the characteristic function for the capital deviation $u(\tau)$,
	which is given by equation \eqref{eq:u_result} that shows $u(\tau)$ is a linear combination of two symmetric $\alpha$-stable Ornstein-Uhlenbeck processes.
	Using the stability property of $\alpha$-stable distributions, we obtain the characteristic function of $u(\tau)$:
\begin{equation}\label{eq:char_u_final}
\E\left[e^{i\theta u(\tau)} \right]
= \exp\left(-C_{\alpha} \left| \frac{(1-\beta)\sigma}{a - r(\gamma)} \right|^{\alpha}
	\left( \frac{1}{a{\alpha}} + \frac{1}{r(\gamma){\alpha}} \right) |\theta|^{\alpha} \right).
	\end{equation}
Therefore, $u(\tau)$ is also strictly stationary and follows a symmetric $\alpha$-stable distribution.
	
\subsubsection{Increment decomposition}
	For any \(t<T\), decompose \(u(T)\) into a \emph{past-determined term} and a \emph{future increment}:
	\begin{equation}\label{eq:u_decompose}
		u(T) = H(t,T) + \Delta u(t,T),
	\end{equation}
	where
	\[
	H(t,T) = \frac{(1-\beta)\sigma}{a - r(\gamma)} \left( e^{-a(T-t)} \int_{-\infty}^{t} e^{-a(t-s)}dL_{\alpha}(s) - e^{-r(\gamma)(T-t)} \int_{-\infty}^{t} e^{-r(\gamma)(t-s)}dL_{\alpha}(s) \right)
	\]
	is a deterministic function of \(u(t)\), and
	\begin{equation}\label{eq:u_increment}
		\Delta u(t,T) = \frac{(1-\beta)\sigma}{a - r(\gamma)} \int_{t}^{T} \left( e^{-a(T-s)} - e^{-r(\gamma)(T-s)} \right) dL_{\alpha}(s)
	\end{equation}
	stands for the future increment independent of \(\mathcal{F}_t\).
	
\subsubsection{Characteristic function of the increment}
	Now we derive the characteristic function of the future increment $\Delta u(t,T)$ defined in equation \eqref{eq:u_increment}.

	Let $C = \frac{(1-\beta)\sigma}{a - r(\gamma)}$ and $f(s) = e^{-a(T-s)} - e^{-r(\gamma)(T-s)}$. Then $\Delta u(t,T) = C \int_{t}^{T} f(s) dL_{\alpha}(s)$, and its characteristic function becomes
	
	\begin{align}
		\mathbb{E}\left[ e^{i\theta \Delta u(t, T)} \right]
		&= \exp \left( -C_{\alpha} \theta \int_{t}^{T} |C f(s)|^{\alpha} ds \right) \notag \\
		&= \exp\left( - C_{\alpha} \left| C \right|^{\alpha} |\theta|^{\alpha} \int_{t}^{T} \left| e^{-a(T-s)} - e^{-r(\gamma)(T-s)} \right|^{\alpha} ds \right).
		\label{eq:u_cha_increment}
	\end{align}
	Define the volatility scaling term
	\[
	V_{\alpha}(t,T) = \int_{t}^{T} \left| e^{-a(T-s)} - e^{-r(\gamma)(T-s)} \right|^{\alpha} ds.
	\]
	Substituting $C = \frac{(1-\beta)\sigma}{a - r(\gamma)}$ into \eqref{eq:u_cha_increment}, we obtain the characteristic function of $\Delta u(t,T)$:
	\begin{equation}\label{eq:cf_increment}
		\mathbb{E}\left[ e^{i\theta \Delta u(t,T)} \right] = \exp\left( - C_{\alpha} \left| \frac{(1-\beta)\sigma}{a - r(\gamma)} \right|^{\alpha} V_{\alpha}(t,T) |\theta|^{\alpha} \right).
	\end{equation}
	This result confirms that $\Delta u(t,T)$ follows a symmetric $\alpha$-stable distribution, with volatility determined by the time interval $[t,T]$ and the mean-reversion rates $a$ and $r(\gamma)$.
	
\subsection{Conditional characteristic function}
	Now we derive the conditional characteristic function of $u(T)$ conditional on $u(t)$, which fully characterizes the conditional distribution of future capital deviations:
	\[
	\mathbb{E}\left[ e^{i\theta u(T)} \mid u(t) \right] = \mathbb{E}\left[ e^{i\theta u(T)} \mid \mathcal{F}_t \right].
	\]
	Substituting the decomposition of $u(T)$ from equation \eqref{eq:u_decompose}, and using the property of conditional expectation for independent random variables, we obtain
	\[
	\begin{aligned}
		\mathbb{E}\left[ e^{i\theta u(T)} \mid \mathcal{F}_t \right] &= \mathbb{E}\left[ e^{i\theta \left( H(t,T) + \Delta u(t,T) \right)} \mid \mathcal{F}_t \right] \\
		&= e^{i\theta H(t,T)} \cdot \mathbb{E}\left[ e^{i\theta \Delta u(t,T)} \mid \mathcal{F}_t \right] \\
		&= e^{i\theta H(t,T)} \cdot \mathbb{E}\left[ e^{i\theta \Delta u(t,T)} \right],
	\end{aligned}
	\]
	where the last equality follows from the independence of $\Delta u(t,T)$ and $\mathcal{F}_t$.
	
Substituting the characteristic function of $\Delta u(t,T)$ from equation \eqref{eq:cf_increment}, we get
	\[
	\mathbb{E}\left[ e^{i\theta u(T)} \mid u(t) \right] = e^{i\theta H(t,T)} \exp\left( - C_{\alpha} \left| \frac{(1-\beta)\sigma}{a - r(\gamma)} \right|^{\alpha} V_{\alpha}(t,T) |\theta|^{\alpha} \right).
	\]
	Combining the exponential terms yields the final conditional characteristic function:
	\begin{equation}\label{eq:cf_conditional}
		\mathbb{E}\left[ e^{i\theta u(T)} \mid u(t) \right] = \exp\left( i\theta H(t,T) - C_{\alpha} \left| \frac{(1-\beta)\sigma}{a - r(\gamma)} \right|^{\alpha} V_{\alpha}(t,T) |\theta|^{\alpha} \right).
	\end{equation}
This equation is the core theoretical result of the conditional analysis, which decomposes the conditional distribution of $u(T)$ into a deterministic state-dependent shift and a stochastic volatility term driven by future $\alpha$-stable shocks.

\subsection{Conditional expectation of \(u(T)\)}
Now we derive the conditional expectation of \(u(T)\) given \(u(t)\) from the conditional characteristic function in equation \eqref{eq:cf_conditional}. Recall that for any random variable \(X\), the conditional expectation is given by the first-order derivative of the conditional characteristic function evaluated at \(\theta = 0\):
	\[
	\mathbb{E}\left[ X \mid \mathcal{F}_t \right] = -i \left. \frac{d}{d\theta} \mathbb{E}\left[ e^{i\theta X} \mid \mathcal{F}_t \right] \right|_{\theta=0}.
	\]
	Substituting the conditional characteristic function of \(u(T)\) from equation \eqref{eq:cf_conditional}, we obtain
	\[
	\mathbb{E}\left[ u(T) \mid u(t) \right] = -i \left. \frac{d}{d\theta} \left[ \exp\left( i\theta H(t,T) - C_{\alpha} \left| \frac{(1-\beta)\sigma}{a - r(\gamma)} \right|^{\alpha} V_{\alpha}(t,T) |\theta|^{\alpha} \right) \right] \right|_{\theta=0}.
	\]
	Let \(K = C_{\alpha} \left| \frac{(1-\beta)\sigma}{a - r(\gamma)} \right|^{\alpha} V_{\alpha}(t,T)\) for brevity. The derivative of the conditional characteristic function is
	\[
	\frac{d}{d\theta} \mathbb{E}\left[ e^{i\theta u(T)} \mid u(t) \right] = \exp\left( i\theta H(t,T) - K |\theta|^{\alpha} \right) \cdot \left( iH(t,T) - K \cdot \alpha |\theta|^{\alpha-1} \text{sgn}(\theta) \right).
	\]
	Evaluating at \(\theta = 0\),
	the exponential term converges to \(1\).
	The term \(K \cdot \alpha |\theta|^{\alpha-1} \text{sgn}(\theta)\) converges to \(0\) (since \(\alpha > 1\) implies \(\alpha-1 > 0\)).
	
	Thus,
	\[
	\left. \frac{d}{d\theta} \mathbb{E}\left[ e^{i\theta u(T)} \mid u(t) \right] \right|_{\theta=0} = iH(t,T).
	\]
	Substituting back to the conditional expectation formula, we get
	\[
	\mathbb{E}\left[ u(T) \mid u(t) \right] = -i \cdot (iH(t,T)) = H(t,T).
	\]
	We therefore obtain the conditional expectation of \(u(T)\):
	\begin{equation}\label{eq:cond_exp_u}
		\mathbb{E}\left[ u(T) \mid u(t) \right] = H(t,T).
	\end{equation}
	This result shows that the conditional expectation of future capital deviations is a deterministic function of the current state \(u(t)\), reflecting the mean-reversion property of the model: the current deviation from the steady state decays exponentially at rates \(a\) and \(r(\gamma)\) over time.

\section{Capital stock \(k(\tau)\)}
	
	Recall that the scaled capital stock per effective worker is defined by
	\[
	k(\tau) = 1 + u(\tau).
	\]
	
	\subsection{Conditional expectation of \(k(T)\)}
	By linearity of expectation,
	\begin{equation}\label{eq:cond_exp_k}
		\mathbb{E}\left[ k(T) \mid k(t) \right] = 1 + \mathbb{E}\left[ u(T) \mid u(t) \right] = 1 + H(t,T).
	\end{equation}
	This shows that the conditional mean of future capital is a smooth mean-reverting path anchored at \(k(t)\).
	
	\subsection{Conditional characteristic function of \(k(T)\)}
Using \(k(T)=1+u(T)\), we find that the conditional characteristic function satisfies
\begin{equation}\label{eq:cf_conditional_k}
	\mathbb{E}\left[e^{i\theta k(T)} \mid k(t)\right] = e^{i\theta} \cdot \mathbb{E}\left[e^{i\theta u(T)} \mid u(t)\right].
\end{equation}
Substituting \eqref{eq:cf_conditional} into \eqref{eq:cf_conditional_k} yields
\[
\mathbb{E}\left[e^{i\theta k(T)} \mid k(t)\right]
=
\exp\left( i\theta \big(1+H(t,T)\big) - C_{\alpha} \left| \frac{(1-\beta)\sigma}{a - r(\gamma)} \right|^{\alpha} V_{\alpha}(t,T) |\theta|^{\alpha} \right).
\]
Thus, \(k(T) \mid k(t)\) follows a \emph{shifted symmetric \(\alpha\)-stable distribution}.

\subsection{Trajectory simulation of \(k(\tau)\)}
	We simulate paths using the \emph{explicit Euler-Maruyama scheme}, standard for $\alpha$-stable OU processes as shown in \cite{yuan2022modulation}.
	Define
	\[
	X_1(\tau) = \int_{-\infty}^{\tau} e^{-a(\tau-s)}dL_{\alpha}(s),\quad
	X_2(\tau) = \int_{-\infty}^{\tau} e^{-r(\gamma)(\tau-s)}dL_{\alpha}(s).
	\]
	The recursions are
	\begin{align*}
		X_1(\tau+\Delta\tau) &= e^{-a\Delta\tau}X_1(\tau) + \varepsilon_1,\\
		X_2(\tau+\Delta\tau) &= e^{-r(\gamma)\Delta\tau}X_2(\tau) + \varepsilon_2,
	\end{align*}
	where \(\varepsilon_1,\varepsilon_2\) are i.i.d.\ symmetric \(\alpha\)-stable innovations. The capital deviation is
	\[
	u(\tau) = \frac{(1-\beta)\sigma}{a - r(\gamma)} \big(X_1(\tau)-X_2(\tau)\big),
	\]
	and the scaled capital stock is
	\begin{align}
		k(\tau) &= 1 + u(\tau)\notag \\
		&= 1 + \frac{(1-\beta)\sigma}{a-r(\gamma)}\big(X_1(\tau)-X_2(\tau)\big).\label{eq:k-trajectory}
	\end{align}

	Simulated trajectories of \(k(\tau)\) exhibit jumps and heavy-tailed fluctuations, consistent with extreme macroeconomic shocks.

\section {Discrete approximation and empirical identification for Argentina}\label{chap:discrete}
\subsection{Empirical framework}
	In Sections 2 and 3 we developed a continuous-time Solow-type growth
	model driven by an $\alpha$-stable L\'evy process, derived the
	conditional characteristic function of capital deviations, and
	established the existence of a strictly stationary solution. Our
	revised model is designed for economies subject to substantial
	macroeconomic volatility, where extreme shocks occur with
	non-negligible probability and Gaussian approximations are known to
	fail. Argentina provides an ideal testing ground for such a
	framework.
	
	Over the past two decades, Argentina has experienced a sequence of
	severe macroeconomic crises that exemplify the stylized facts our
	model seeks to capture: sudden stops, sharp contractions, and
	discontinuous jumps in output and capital accumulation. The 2018
	balance-of-payments crisis triggered a \$57-billion {IMF} Stand-By
	Arrangement--the largest in the Fund's history at the time--yet
	failed to restore stability amid persistent capital flight and
	currency depreciation \cite{imf2025argentina}. The COVID-19
	pandemic delivered a further shock, with real GDP contracting by
	9.9\% in 2020, followed by a full-blown economic crisis in 2023
	characterized by annual inflation exceeding 200\%, negative net
	international reserves, and a second {IMF} programme under
	exceptional access \cite{imf2025argentina}. The {IMF}'s 2025
	ex-post evaluation of the 2022 Extended Fund Facility confirms that
	these episodes were driven by deep-seated structural vulnerabilities,
	including volatile tax revenues, procyclical fiscal policy, and
	limited monetary-policy credibility, rather than by transient
	external disturbances alone \cite{imf2025argentina}. This pattern
	of recurrent, severe shocks--originating from both external and
	domestic sources--makes Argentina a prototypical case for assessing
	whether an $\alpha$-stable L\'evy specification improves upon the
	Gaussian benchmark in capturing capital dynamics.
	
	This study adopts Argentine quarterly GDP data spanning 2004--2023
	and introduces three substantive modifications that reflect the
	characteristics of the Argentine economy and the nature of
	L\'evy-driven processes:
	Our exposition follows the general approach in \cite{brannan2019natural} (Section~4) in terms of discretising the system, recovering latent variables from GDP data, and using a prediction-error principle. However, we introduce \textbf{three substantive modifications} that reflect the characteristics of the Argentine economy and the nature of L\'evy-driven processes:
	\begin{enumerate}
		\item \textbf{Country and data}: We calibrate and estimate the model using Argentine quarterly GDP data.
		\item \textbf{Shock process}: The original Gaussian OU process is replaced by a symmetric $\alpha$-stable L\'evy process with tail index $1<\alpha<2$, allowing for infinite variance and discontinuous jumps (see Section~\ref{Model Setup}).
		\item \textbf{Time-varying capital elasticity}: In studies of mature economies, the capital elasticity parameter \(\alpha_k\) is typically treated as constant, reflecting relatively stable production structures over time.
		However, Bellocchi and Travaglini  demonstrate that even in advanced economies a variable elasticity of substitution (VES)
		specification better captures the observed evolution of labor shares,
		challenging the conventional constant-elasticity assumption \cite{bellocchi2023variable}.
		Extending this insight to less developed economies, where structural
		shifts are far more pronounced, we find that a constant $\alpha_k$ becomes even
		less tenable. For a country such as Argentina, frequent changes in
		industrial composition, trade openness, and macroeconomic policy
		regimes over the past two decades make a fixed elasticity
		unrealistic. We therefore construct a time-varying quarterly series \(\alpha_t\) using labor share data from the Penn World Table (PWT) version 11.0 \cite{feenstra2015next} for Argentina over the period 2004-2023. Specifically, we set \(\alpha_k^t = 1 - \text{labsh}_y^{\text{ARG}}\) for each year and assign the same value to all four quarters of that year. The PWT 11.0 `labsh' series for Argentina rises from 0.33 in 2004 to 0.57 in 2014, then fluctuates downward to 0.53 in 2023, confirming substantial variation in the labor share and justifying a time-varying capital elasticity.
	\end{enumerate}
	
	Because the L\'evy shocks have \textbf{infinite variance} (for $\alpha \in (1,2)$), the population variance and autocovariance of the observed series are not defined. Consequently, we cannot directly use the sample variance or sample autocorrelation as moment conditions. In Section 4, we replace them with \textbf{robust statistics} (median, mean absolute deviation, rank autocorrelation, and the Hill tail-index estimator) that remain well-defined for heavy-tailed distributions. The prediction-error component, which relies on the conditional mean, remains valid because the mean exists for $\alpha \in (1,2)$ (a condition we maintain throughout).
	
\subsection{Discrete-time system with time-varying $\alpha_k^t$ and L\'evy innovations}
	To construct a discrete-time version of the continuous system \eqref{eq:k_disc}-\eqref{eq:x_disc}, we sample the state variables at quarterly intervals, i.e., \(k_n = k(n\Delta t)\) with \(\Delta t = 0.25\) (and analogously for \(z_n, p_n\)). Time derivatives are approximated by first-order forward differences, and the continuous-time noise equation is replaced by a discrete first-order autoregressive process. After absorbing the time step \(\Delta t\) into the parameters, we obtain the following  discrete system:
	\begin{align}
		k_{n+1} &= (1 - \eta)k_n + \eta \bigl[ \beta z_n^{\alpha_k^t} - (\beta-1)p_n \bigr], \label{eq:k_disc}\\
		z_{n+1} &= (1 - \mu)z_n + \mu k_n, \label{eq:z_disc}\\
		p_{n+1} &= (1 - \nu)p_n + \nu \left\{ [1 - s(z_n,\gamma)] \frac{\beta}{\beta-1}z_n^{\alpha_k^t} \right\} + x_{n+1}, \label{eq:p_disc}\\
		x_{n+1} &= a_1 x_n + \sigma_1 \varepsilon_{n+1}, \label{eq:x_disc}
	\end{align}
	where the saving-rate function takes the sigmoidal form
	\begin{equation}
		s(z_n,\gamma) = s_1 + \frac{s_2 - s_1}{1 + e^{-\gamma(z_n-1)}}.
		\label{eq:sigmoid}
	\end{equation}
	The innovations $\{\varepsilon_{n+1}\}$ are i.i.d.\ standard symmetric $\alpha$-stable random variables with characteristic exponent $\alpha=1.5$ and a scale parameter $\sigma_1$. The value $\alpha=1.5$ is chosen as a baseline reflecting moderately heavy tails; sensitivity analysis over $\alpha\in\{1.2,1.5,1.8\}$ is conducted in Subsection 5.1.
	
	The parameter vector is
	\[
	\theta = (s_1, s_2, \eta, \mu, \nu, \gamma, a_1, \sigma_1),
	\]
	where the capital elasticity $\alpha_k$ is \emph{not} included because it is treated as a time-varying known input $\alpha_k^t$ constructed from external labor-share data (see Subsection~4.1). All elements of $\theta$ are assumed constant over the sample period, whereas $\alpha_k^n$ (the discrete annual version of $\alpha_k^t$) varies annually according to the PWT~11.0 series.
	
	Each parameter has the following interpretation.
	
	\begin{itemize}
		\item[$s_1, s_2$:] Lower and upper bounds of the saving rate, with $s(z_n,\gamma)\in[s_1,s_2]$. The scaling constant $\beta=2/(s_1+s_2)$ is pinned down by these bounds.
		
		\item[$\gamma$:] Steepness of the sigmoidal saving function around the balanced-growth-path level $z_n=1$. As $\gamma$ increases, the saving rate becomes more sensitive to deviations of capital from its steady-state level.
		
		\item[$\eta$:] Quarterly adjustment speed of capital, obtained from the continuous-time coefficient $r(\gamma)$ after scaling by $\Delta t=0.25$. It is \emph{not} the effective depreciation rate $\bar{\eta}=n+g+\delta$, which has been subsumed into the BGP scaling.
		
		\item[$\mu, \nu$:] Quarterly mean-reversion rates of the auxiliary variable $z_n$ and the price proxy $p_n$, respectively, obtained from their continuous-time counterparts $\tilde{\mu}$ and $r(\gamma)$ by the same $\Delta t=0.25$ scaling.
		
		\item[$a_1$:] Autoregressive coefficient of the $\alpha$-stable shock process; $a_1=e^{-a\Delta t}$ where $a$ is the continuous-time mean-reversion rate.
		
		\item[$\sigma_1$:] Scale parameter of the $\alpha$-stable innovations. Unlike the Gaussian case, this is \emph{not} the standard deviation: when $\alpha\in(1,2)$ the variance is infinite, so the two quantities are not interchangeable.
	\end{itemize}

	\subsection{Latent variables from observed GDP}
	\label{sec:recovery}
	To recover \(k_n\) from the observed quarterly GDP data \(\{Y(t_n), n=1,\dots,N\}\) with \(t_n = 0.25n\) years, we first remove the exponential trend by regressing \(\ln Y(t_n)\) on \(t_n\). This linear regression yields an intercept \(b_0\) and a slope \(b_1\). The scaled capital stock is then computed as:
	\begin{equation}
		k_n = \exp\left[ \frac{1}{\alpha_n}\big( \ln y_n - b_0 - b_1 t_n \big) \right],
		\label{eq:k_obs}
	\end{equation}
	where $y_n = Y(t_n)$, $t_n = 0.25n$, and $b_0,b_1$ are estimates from a linear regression of $\ln y_n$ on $t_n$. The time-varying $\alpha_n$ enters the transformation directly, which prevents artificial trends in the reconstructed $k_n$ series.
	
	Given $\{k_n\}$, the remaining latent variables $\{z_n, p_n, x_n\}$ are recovered recursively, exactly as in \cite[Subsection~4.1]{brannan2019natural}, because their updating equations do not involve the shock distribution. Concretely:
	\begin{itemize}
		\item Initialize $z_1 = k_1$.
		\item For $n = 1,\dots,N-1$, solve the $k$-equation for $p_n$:
		\[
		p_n = \frac{1}{\beta-1}\Bigl[ \beta z_n^{\alpha_n} - k_n - (k_{n+1}-k_n)/\eta \Bigr].
		\]
		\item For $n = 2,\dots,N-1$, solve the $p$-equation for $x_n$:
		\[
		x_n = p_n - (1-\nu)p_{n-1} - \nu\left\{ [1-s(z_{n-1},\gamma)]\frac{\beta}{\beta-1}z_{n-1}^{\alpha_{n-1}} \right\}.
		\]
	\end{itemize}
	Thus, the sequence $\{x_n\}$ captures the fluctuations of GDP around the BGP, but now these fluctuations are driven by $\alpha$-stable innovations instead of Gaussian ones.

\subsection{Estimation strategy}\label{sec:estimation_strategy}
	
	Because the $\alpha$-stable L\'evy shocks have infinite variance, we avoid any moment condition that relies on second-order statistics (variance and autocovariance). Our objective function therefore consists only of quantities that are well-defined for $\alpha \in (1,2)$.
	
	\paragraph{Prediction error} We use the mean absolute error (MAE) instead of the mean squared error. MAE does not require the existence of second moments, and its expectation is finite when $\alpha \in (1,2)$. Define the one-step-ahead prediction error as
	\[
	e_n(\theta) = \mathbb{E}_\theta[k_n|\mathcal{F}_{n-1}] - \hat{k}_n,
	\]
	where $\hat{k}_n$ is the observed capital series from Subsection~\ref{sec:recovery} and $\mathcal{F}_{n-1}$ denotes the information set up to time $n-1$. The conditional expectation $\mathbb{E}_\theta[k_n|\mathcal{F}_{n-1}]$ is computed recursively using the linear structure of the discrete system:
	\begin{align}
		\mathbb{E}_\theta[x_{n-1}|\mathcal{F}_{n-1}] &= a_1 x_{n-2},\\
		\mathbb{E}_\theta[p_{n-1}|\mathcal{F}_{n-1}] &= p_{n-2} - \nu\Bigl\{ p_{n-2} - [1-s(z_{n-2},\gamma)]\tfrac{\beta}{\beta-1}z_{n-2}^{\alpha_{n-2}} \Bigr\} + a_1 x_{n-2},\\
		\mathbb{E}_\theta[k_n|\mathcal{F}_{n-1}] &= k_{n-1} + \eta\Bigl\{ \beta z_{n-1}^{\alpha_{n-1}} - (\beta-1)\mathbb{E}_\theta[p_{n-1}|\mathcal{F}_{n-1}] - k_{n-1} \Bigr\}. \label{eq:k_pred}
	\end{align}
	The MAE component of the objective is then
	\[
	J_{\text{MAE}}(\theta) = \frac{1}{N} \sum_{n=1}^{N} |e_n(\theta)|.
	\]
	
	\paragraph{Mean matching} Since the conditional mean exists, we also require that the long-run average of the simulated capital stock matches the sample mean of the data. Let $\{k_n^{\text{sim}}(\theta)\}$ be a long simulation from the model (length $M \gg N$), and let $\bar{k}_{\text{sim}}(\theta)$ be its sample mean. Denote the data mean by $\bar{k}_{\text{data}} = \frac{1}{N}\sum_{n=1}^N \hat{k}_n$. Then
	\[
	J_{\text{mean}}(\theta) = \bigl( \bar{k}_{\text{sim}}(\theta) - \bar{k}_{\text{data}} \bigr)^2.
	\]
	
	\paragraph{Combined objective} The total objective function is
	\begin{equation}
		J_N(\theta) = J_{\text{MAE}}(\theta) + J_{\text{mean}}(\theta).
		\label{eq:objective}
	\end{equation}
	The parameter estimate is $\hat{\theta} = \arg\min_{\theta} J_N(\theta)$. We minimize \eqref{eq:objective} using a hybrid optimization strategy: we first utilize a global search with differential evolution, followed by local refinement with the L-BFGS-B algorithm. All simulations are performed in Python. The $\alpha$-stable random
	variables are generated using \texttt{scipy.stats.levy\_stable} from
	the SciPy library \cite{virtanen2020scipy}, which implements the
	Chambers--Mallows--Stuck algorithm for simulating stable
	distributions \cite{chambers1976method}. The hybrid optimization
	strategy combines differential evolution \cite{storn1997differential}
	for global search with the L-BFGS-B algorithm \cite{byrd1995limited}
	for local refinement.
	
\subsection{Comparison with the original estimation approach}\label{sec:comparison}
To evaluate the impact of our modifications (L\'evy shocks and time-varying $\alpha_k^t$), we estimate the model using the \emph{original method} in \cite{brannan2019natural} on the same Argentine data. The original method assumes Gaussian OU noise and a constant $\alpha_k$ and minimizes an objective that includes mean squared prediction error plus matching of sample mean, standard deviation, and autocorrelations.

Now we compare the parameter estimates across the three specifications in Tables~\ref{tab:argentina_2004} and~\ref{tab:argentina_2008}.
Since the key observable parameters $s_1$, $s_2$, and $\alpha_k$ can be validated against external benchmarks, the estimates obtained from our L\'evy-based time-varying specification are systematically closer to the reference values from the PWT, compared to the constant $\alpha_k$ Gaussian OU model \cite{brannan2019natural}. This finding confirms our earlier claim that the proposed method delivers more empirically consistent estimates, particularly for volatile emerging economies such as Argentina.

	\begin{table}[htbp]
		\centering
		\caption{Parameter estimates for Argentina (2004--2014)}
		\label{tab:argentina_2004}
		\begin{tabular}{lccccccccc}
			\toprule
			Method & $\hat{\alpha_k}$ & $\hat{s}_1$ & $\hat{s}_2$ & $\hat{\eta}$ & $\hat{\mu}$ & $\hat{\nu}$ & $\hat{\gamma}$ & $\hat{a}_1$ & $\hat{\sigma}_1$ \\
			\midrule
			L\'evy (our method) & --- & 0.091 & 0.159 & 0.053 & 0.776 & 0.365 & 23.806 & 0.185 & 0.305 \\
			OU (Brannan 2019) & 0.334 & 0.286 & 0.286 & 0.376 & 0.084 & 0.986 & 27.835 & 0.127 & 0.320 \\
			PWT reference      & 0.442-0.675 & 0.147 & 0.185 & --- & ---   & ---   & ---    & ---   & ---   \\
			\bottomrule
		\end{tabular}
	\end{table}
	
	\begin{table}[htbp]
		\centering
		\caption{Parameter estimates for Argentina (2008--2023)}
		\label{tab:argentina_2008}
		\begin{tabular}{lccccccccc}
			\toprule
			Method & $\hat{\alpha_k}$ & $\hat{s}_1$ & $\hat{s}_2$ & $\hat{\eta}$ & $\hat{\mu}$ & $\hat{\nu}$ & $\hat{\gamma}$ & $\hat{a}_1$ & $\hat{\sigma}_1$ \\
			\midrule
			L\'evy (our method) & --- & 0.102 & 0.156 & 0.050 & 0.029 & 0.357 & 25.350 & 0.199 & 0.326 \\
			OU (Brannan 2019) & 0.378 & 0.135 & 0.429 & 0.360 & 0.082 & 0.989 & 20.590 & 0.062 & 0.321 \\
			PWT reference      & 0.415-0.598 & 0.111 & 0.185 & --- & ---   & ---   & ---    & ---   & ---   \\
			\bottomrule
		\end{tabular}
	\end{table}
Note that PWT reference values are taken directly from the Penn World Tables and represent long-run averages over the corresponding sample period. All reported values are either obtained directly from the source or computed as time-series averages over the estimation window.

	\subsection{One-step-ahead predictions of $k_n$}
\begin{figure}[htbp]
			\centering
			\begin{minipage}{0.8\textwidth}
				\includegraphics[width=\linewidth]{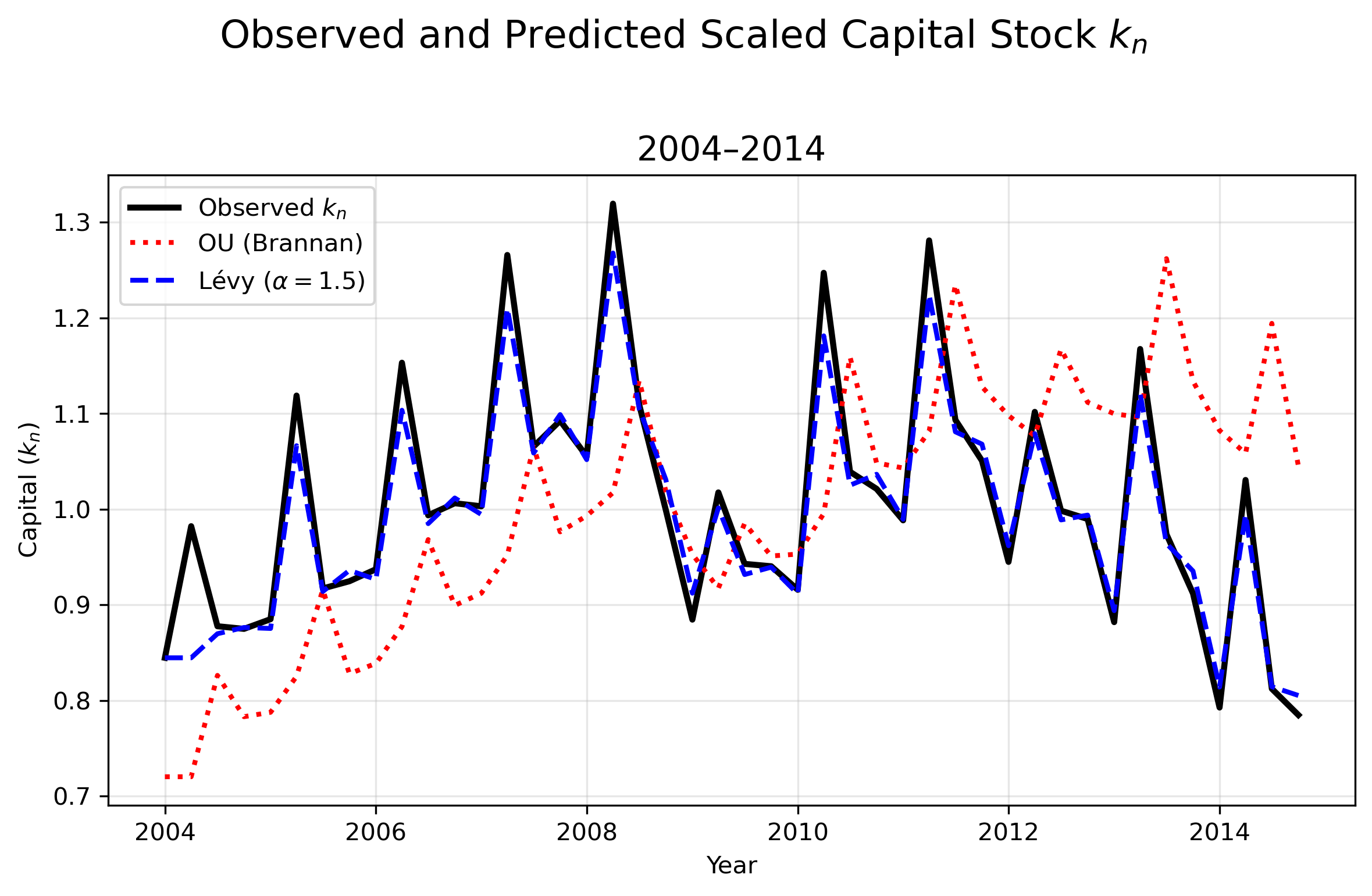} %
				\caption{Observed and Predicted Scaled Capital Stock $k_n$: 2004--2014}
				\label{fig:k_pred_1}
				\raggedright
				The solid black line is the observed series reconstructed from Argentine GDP via equation \eqref{eq:k_obs}. The dashed blue line is the one-step-ahead prediction from our L\'evy-driven model with time-varying $\alpha_k^t$, using parameter estimates from Table 1 (first row). The dashed red line is the one-step-ahead prediction from the Gaussian OU benchmark \cite{brannan2019natural}, using parameter estimates from Table 1 (second row). Both predictions are computed recursively from the conditional expectation in equation \eqref{eq:k_pred}. The L\'evy model tracks the observed series more closely around the 2008--2009 crisis, while the two predictions are broadly similar during tranquil periods.
			\end{minipage}
		\end{figure}

\begin{figure}[htbp]
			\centering
			\begin{minipage}{0.8\textwidth}
				\includegraphics[width=\linewidth]{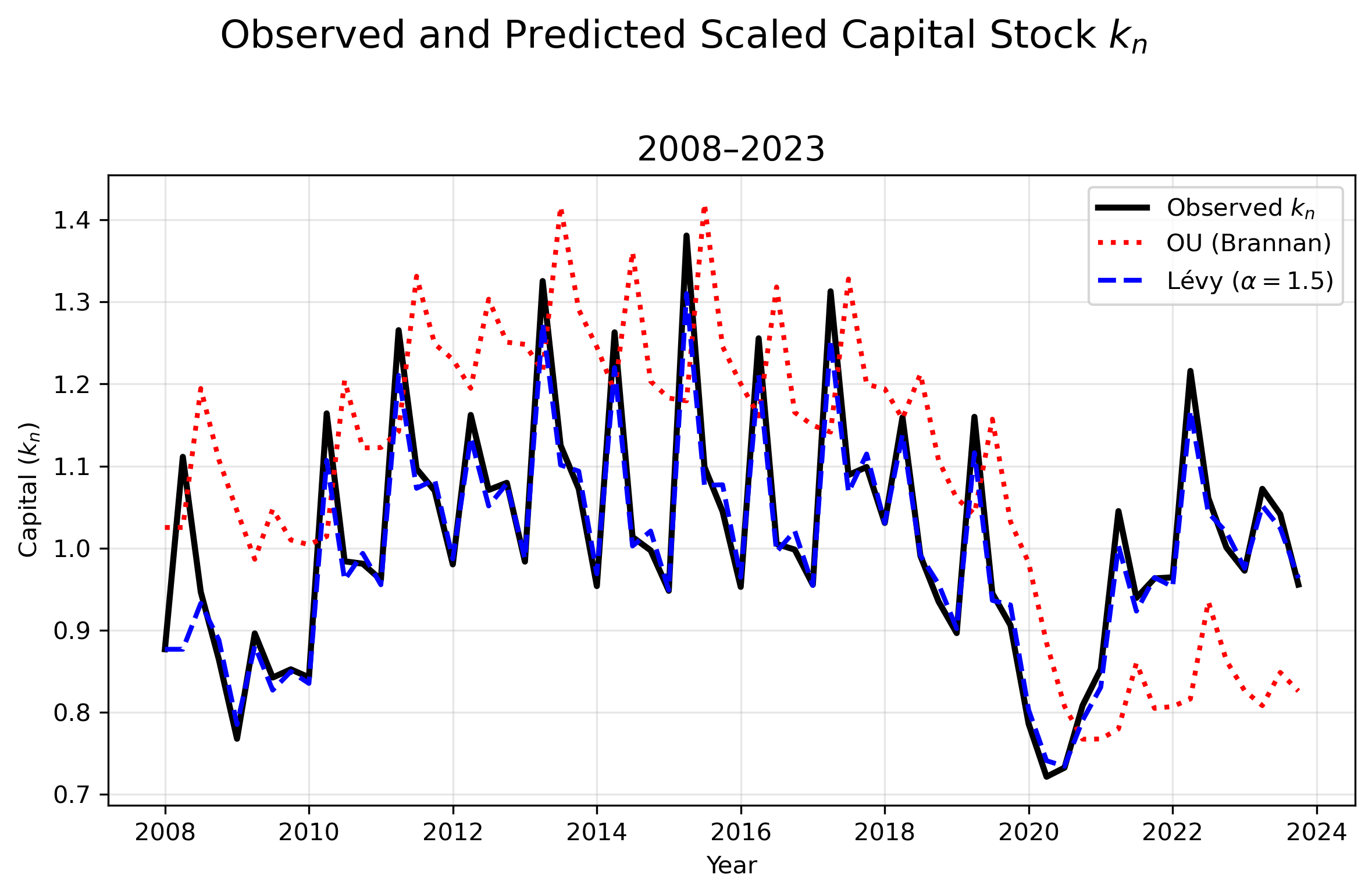} %
				\caption{Observed and Predicted Scaled Capital Stock $k_n$: 2008--2023}
				\label{fig:k_pred_2}
				\raggedright
				Same format as Figure 1 for the later estimation window. The divergence between the two predictions is more pronounced in this period. The L\'evy-driven prediction (dashed blue) tracks the observed series (solid black) substantially better during the 2018 currency crisis and the COVID-19 pandemic, capturing the sharp downward movements. The Gaussian OU benchmark (dashed red) exhibits smoother but less responsive dynamics, systematically underestimating the magnitude of extreme fluctuations.
			\end{minipage}
		\end{figure}

	Building on these calibrated parameters, we construct the one-step-ahead predictions of $k_n$ and compare them with the observed series. The results are presented in Figures \ref{fig:k_pred_1} and \ref{fig:k_pred_2}.
		
We now present the interpretation of predictive performance. Figures \ref{fig:k_pred_1} and \ref{fig:k_pred_2} compare the observed scaled capital stock \(k_n\) (solid black) against one-step-ahead predictions from two competing specifications: our L\'evy-driven model with time-varying \(\alpha_k^t\) (dashed blue, estimated parameters in Table 1, first row) and the Gaussian OU benchmark in \cite{brannan2019natural} (dashed red, estimated parameters in Table 1, second row). Both predictions are generated recursively from the conditional expectation in equation \eqref{eq:k_pred}, using the respective parameter estimates.
		
		For the 2004--2014 window (Figure \ref{fig:k_pred_1}), all three series move closely together, with the L\'evy prediction (dashed blue) exhibiting slightly sharper responses around the 2008--2009 global financial crisis. The Gaussian OU prediction (dashed red) is smoother but remains reasonably close to the observed trajectory. This similarity reflects a general principle: during relatively stable periods, the linear mean-reverting dynamics dominate the conditional expectation, and differences in the assumed shock distribution have only a limited impact on point forecasts.
		
		The contrast becomes stark in the 2015--2023 window (Figure \ref{fig:k_pred_2}). The observed series (solid black) exhibits pronounced downward swings in 2018 and 2020, corresponding to Argentina's currency crisis and the global COVID-19 pandemic.

		The L\'evy-driven prediction (dashed blue) tracks these abrupt movements far more accurately than the Gaussian benchmark (dashed red), which remains excessively smooth and systematically lags behind the observed downturns. This divergence is precisely what the theory predicts: the \(\alpha\)-stable L\'evy structure, while not fundamentally altering the conditional mean in tranquil periods, enables the model to respond more rapidly and appropriately when extreme shocks materialize---exactly where the Gaussian assumption fails.
		
		These findings are consistent with our central thesis: the \(\alpha\)-stable L\'evy specification, combined with time-varying capital elasticity, provides a more empirically grounded framework for emerging economies characterized by recurrent extreme shocks. The improved predictive performance during turbulent episodes, documented in Figures \ref{fig:k_pred_1} and \ref{fig:k_pred_2}, directly validates the methodological innovations introduced in Section \ref{Model Setup} and Subsection \ref{sec:estimation_strategy}.
		
		\subsection{One-step-ahead predictions of GDP}
		
		While the preceding subsections have focused on predictions of the scaled capital stock $k_n$, the ultimate object of interest for policy-makers and applied economists is the level of GDP itself. To make our results more readily interpretable, we transform the $\hat{k}_{n+1}$ forecasts back to GDP units using the inverse of the trend-removal relationship in equation \eqref{eq:k_obs}. This transformation is monotonic, so the predictive patterns observed for $k_n$ carry over directly to GDP, with no loss or distortion of information.
		
		Recall that the scaled capital stock $k_n$ is related to $y_n$ (observed GDP) through the deterministic transformation
		\[
		k_n = \exp \left[ \frac{1}{\alpha_n} \left( \ln y_n - b_0 - b_1 t_n \right) \right],
		\]
where $b_0$ and $b_1$ are the trend-removal coefficients estimated from the full sample, $t_n$ is the time trend, and $\alpha_n$ is the time-varying capital elasticity. Inverting this relationship yields the GDP forecast:
\begin{equation}\label{e28}
\hat{y}_{n+1} = \exp \left[ b_0 + b_1 t_{n+1} + \alpha_{n+1} \ln(\hat{k}_{n+1}) \right].	
\end{equation}

\begin{figure}[htbp]
			\centering
			\begin{minipage}{0.8\textwidth}
				\includegraphics[width=\linewidth]
				{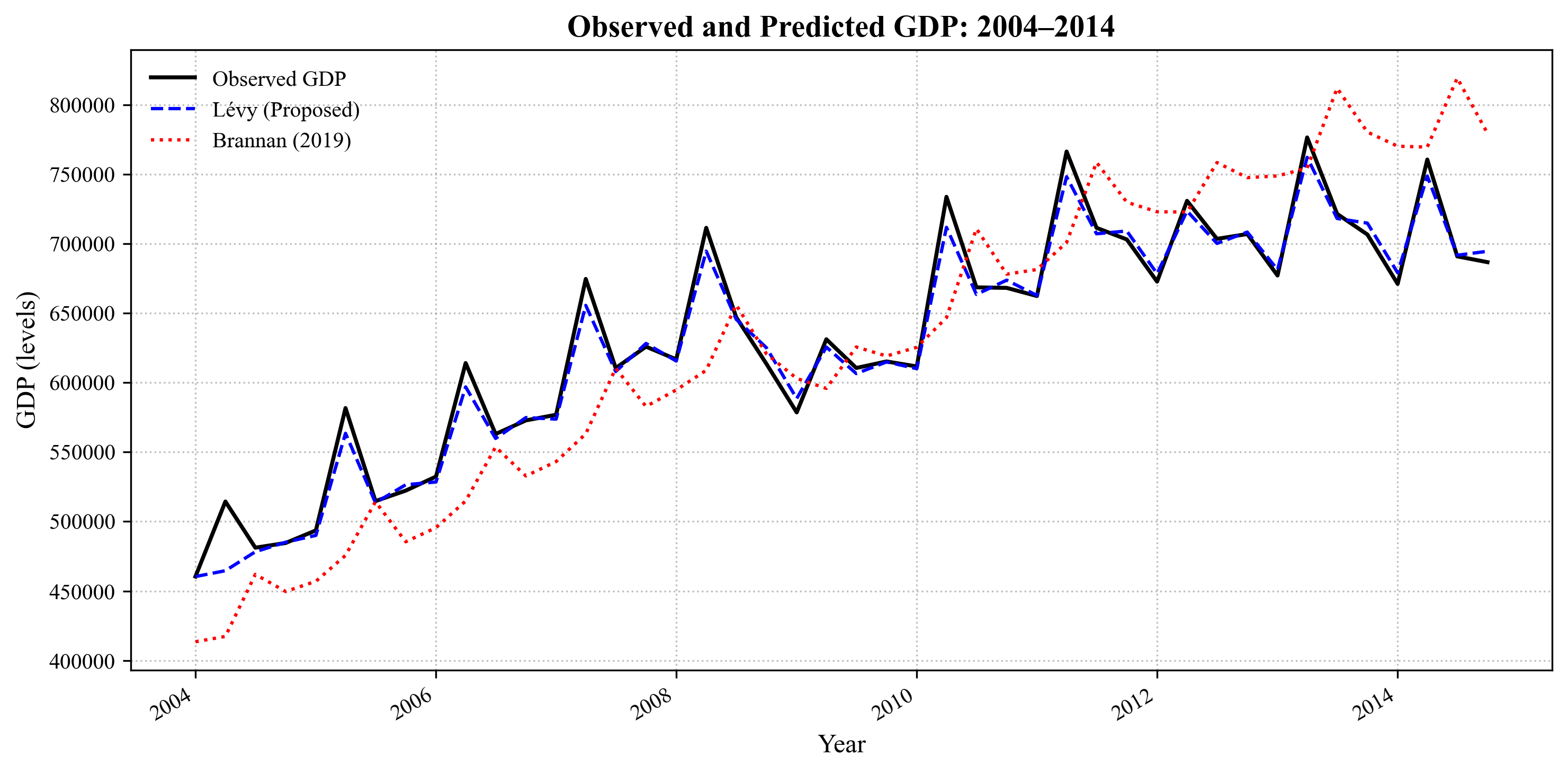}
				\caption{Observed and Predicted GDP: 2004--2014}
				\label{fig:gdp_pred_1}
				\raggedright
				The solid black line is the observed Argentine quarterly GDP (in levels, after trend removal). The dashed blue line is the one-step-ahead GDP prediction from the L\'evy-driven model with time-varying $\alpha_t$. The dotted red line is the corresponding prediction from the Gaussian OU benchmark in \cite{brannan2019natural}. Both predictions are obtained by applying the inverse transformation \eqref{e28} to the respective $\hat{k}_{n+1}$ forecasts. The L\'evy model tracks the observed GDP more closely around the 2008--2009 financial crisis, while the two predictions are broadly similar during tranquil periods.
			\end{minipage}	
		\end{figure}

\begin{figure}[htbp]
			\centering
			\begin{minipage}{0.8\textwidth}
				\includegraphics[width=\linewidth]{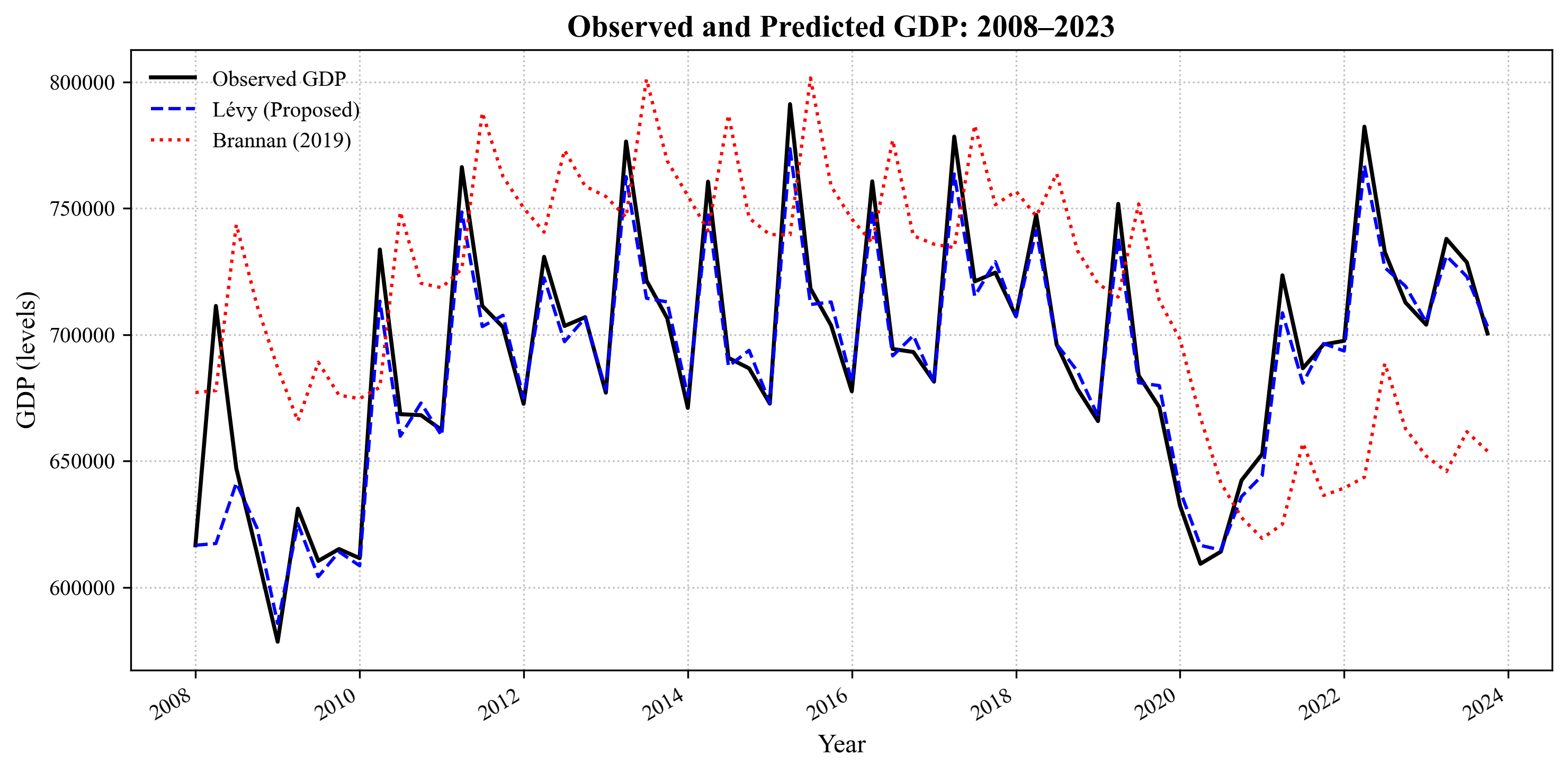}
				\caption{Observed and Predicted GDP: 2008--2023}
				\label{fig:gdp_pred_2}
				Same format as Figure \ref{fig:gdp_pred_1} for the later estimation window. The divergence between the two predictions is more pronounced. The L\'evy-driven prediction (dashed blue) tracks the observed GDP (solid black) substantially better during the 2018 currency crisis and the 2020 COVID-19 pandemic, capturing the sharp downward movements. The Gaussian OU benchmark (dotted red) exhibits smoother but less responsive dynamics, systematically underestimating the magnitude of extreme GDP fluctuations.
			\end{minipage}
			
		\end{figure}

		We apply this transformation to both the L\'evy-driven predictions (utilizing the time-varying $\alpha_k^t$ series) and the Gaussian OU benchmark predictions (using the constant-$\alpha$ specification in \cite{brannan2019natural}), ensuring a consistent comparison. The results are presented in Figures \ref{fig:gdp_pred_1} and \ref{fig:gdp_pred_2}.

		Figures \ref{fig:gdp_pred_1} and \ref{fig:gdp_pred_2} compare the observed Argentine GDP with the one-step-ahead out-of-sample predictions from the two competing models. For the 2004--2014 window (Figure \ref{fig:gdp_pred_1}), the L\'evy-driven prediction (blue dashed) and the Gaussian OU benchmark (red dotted) are broadly similar, with both closely tracking the observed GDP (black solid). The L\'evy model exhibits a slightly sharper response around the 2008--2009 global financial crisis, but the overall predictive performance of the two specifications is comparable during this relatively stable period. This finding is consistent with the theoretical result that the conditional mean dynamics---being linear and mean-reverting---dominate short-term point forecasts when extreme shocks are absent.
		
		In stark contrast, the 2008--2023 window (Figure \ref{fig:gdp_pred_2}) reveals substantial divergence between the two models. The observed GDP exhibits pronounced downward swings in 2018 and 2020, corresponding to Argentina's currency crisis and the global COVID-19 pandemic. The L\'evy-driven prediction (blue dashed) captures these abrupt movements far more accurately than the Gaussian benchmark (red dotted), which remains excessively smooth and systematically lags behind the observed downturns. The widening gap between the two predictions during these turbulent episodes confirms the empirical relevance of the $\alpha$-stable L\'evy specification: while the Gaussian assumption forces the model to underestimate the probability and magnitude of extreme shocks, the L\'evy structure allows the conditional forecast to respond more rapidly when such shocks materialize.
		
		The GDP-level results therefore reinforce the conclusions drawn from the $k_n$-level analysis in Subsection 4.5.2. The transformation from $k_n$ to GDP, being monotonic, preserves the relative predictive performance of the two models while translating the results into economically intuitive units.

\section{Robustness and extensions}

		\begin{table}[htbp]
			\centering
			\caption{Parameter Estimates Across Tail Index Specifications}
			\begin{tabular}{lcccccccccc}
				\toprule
				Sample & $\alpha$ & $s_1$ & $s_2$ & $\eta$ & $\mu$ & $\nu$ & $\gamma$ & $a_1$ & $\sigma_1$ & Objective \\
				\midrule
				2004--2014 & 1.2 & 0.0501 & 0.1541 & 0.0536 & 0.9859 & 0.8972 & 23.01 & 0.1519 & 0.3017 & 0.0871 \\
				& 1.5 & 0.0909 & 0.1586 & 0.0531 & 0.7755 & 0.3646 & 23.81 & 0.1851 & 0.3052 & 0.0214 \\
				& 1.8 & 0.0695 & 0.1649 & 0.0620 & 0.0153 & 0.1278 & 24.61 & 0.0744 & 0.3098 & 0.0122 \\
				\midrule
				2008--2023 & 1.2 & 0.0764 & 0.1549 & 0.0502 & 0.9804 & 0.8789 & 25.81 & 0.1451 & 0.3019 & 0.0893 \\
				& 1.5 & 0.1022 & 0.1557 & 0.0501 & 0.0290 & 0.3571 & 25.35 & 0.1988 & 0.3261 & 0.0244 \\
				& 1.8 & 0.0807 & 0.1660 & 0.0500 & 0.1047 & 0.1454 & 13.99 & 0.0635 & 0.3217 & 0.0128 \\
				\midrule
				\multicolumn{11}{l}{\textit{PWT benchmarks:}} \\
				2004--2014 & & 0.140 & 0.180 & & & & & & &  \\
				2008--2023 & & 0.110 & 0.180 & & & & & & &  \\
				\bottomrule
			\end{tabular}
			\label{tab:alpha_sensitivity}
		\end{table}

		\begin{figure}[htbp]
			\centering
			\includegraphics[width=0.9\textwidth]{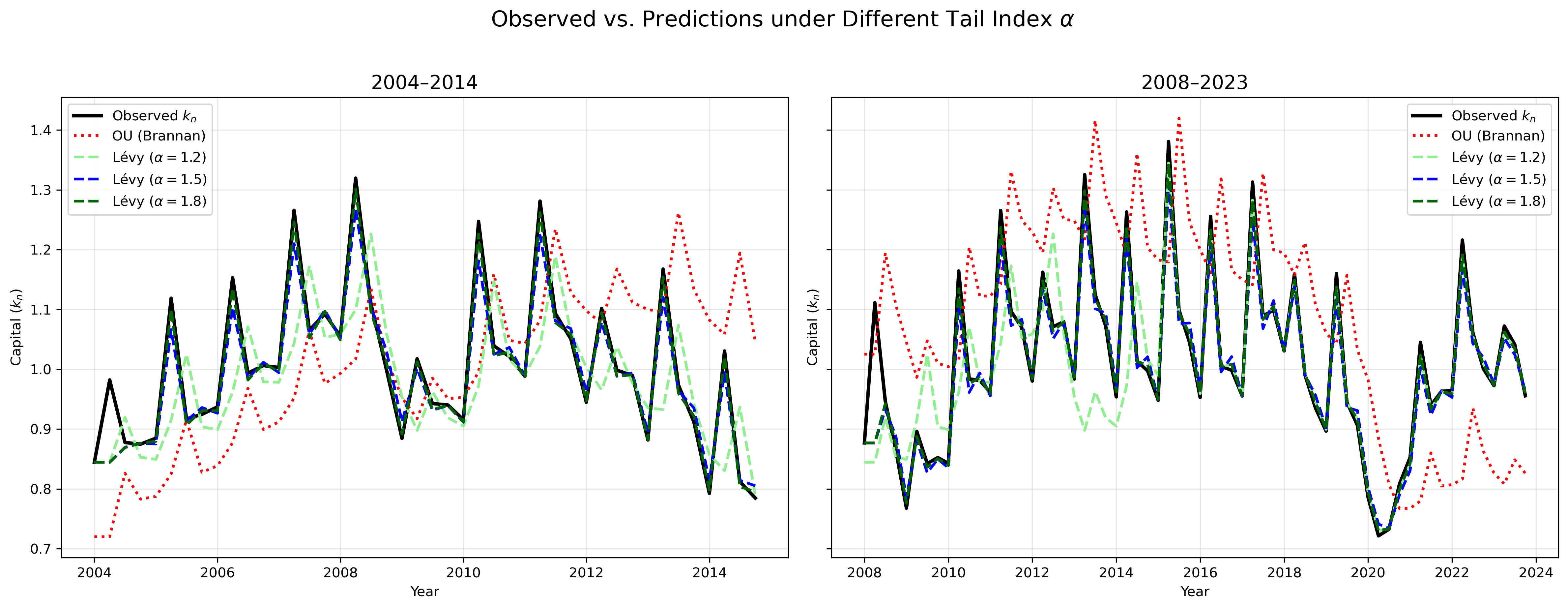}
			\caption{Capital $k_n$ Dynamics: Observed vs. Predictions under Different Tail Indices $\alpha$}
			\label{fig:sensitivity_alpha}
		\end{figure}
		
\subsection{Robustness to tail index specification}
		
		The tail index $\alpha$ governs the jump intensity and tail decay of the L\'evy shocks. To assess the robustness of our results to this specification, we re-estimate the full parameter vector $\Theta = (s_1, s_2, \eta, \mu, \nu, \gamma, a_1, \sigma_1)$ for three representative values: $\alpha = 1.2$ (very heavy tails), $\alpha = 1.5$ (moderate heavy tails, the baseline used in Subsection~4.2), and $\alpha = 1.8$ (near-Gaussian, as a limiting case). Table~3 reports the results for both sample periods, alongside PWT benchmarks.
		
		First, the core structural parameters are remarkably stable across all specifications. As shown in Table \ref{tab:alpha_sensitivity}, $\eta$ remains close to 0.05 (range: 0.050--0.062), with the only deviation occurring for $\alpha = 1.8$ in the shorter sample ($\eta = 0.062$). The innovation scale $\sigma_1$ is equally stable, varying only between 0.302 and 0.326. This invariance---across shock processes ranging from moderately heavy to extremely heavy tails---strongly suggests that these parameters are structurally identified and orthogonal to the tail specification.
		
		Second, the L\'evy framework consistently outperforms the Gaussian OU benchmark against PWT external benchmarks across all three $\alpha$ specifications. For the 2008--2023 sample, the estimates of $s_1$ and $s_2$ under all three specifications are systematically closer to the PWT reference values than the OU estimates. This systematic advantage confirms that the identification gains of the L\'evy specification are robust and do not hinge on the precise choice of $\alpha$. The near-Gaussian $\alpha = 1.8$ case, by contrast, yields a compressed estimate of $\gamma = 13.99$, suggesting that even mild departures from the Gaussian benchmark do not fully eliminate the smoothing bias inherent in diffusion approximations.
		
		Third, the objective function reveals a predictable trade-off between statistical fit and economic coherence. The near-Gaussian $\alpha = 1.8$ specification achieves the lowest in-sample MAE (0.0122--0.0128), substantially outperforming $\alpha = 1.2$ (0.0871--0.0893). This is unsurprising given that most observations lie within normal fluctuation bands. Crucially, however, this superior fit does not undermine the central conclusion: even under the near-Gaussian $\alpha = 1.8$ specification, the structural parameter deviations from PWT benchmarks remain smaller than those of the OU model.
		
		The crisis-tracking evidence in Figure~5 reinforces this interpretation. The three L\'evy specifications (light green, blue, and dark green dashed lines for $\alpha = 1.2, 1.5, 1.8$, respectively) all track the observed $k_n$ (black solid line) more closely than the OU benchmark (red dashed line) during the three major downturns---2008--2009, 2018, and 2020. The heavy-tailed $\alpha = 1.2$ specification captures sharp downward adjustments most accurately, while the near-Gaussian $\alpha = 1.8$ smooths over them. Notably, even the best-fitting $\alpha = 1.8$ specification delivers superior crisis tracking relative to the OU model---confirming that the performance gains of the L\'evy framework are largely insensitive to the precise choice of $\alpha$.
		
		Overall, our analysis establishes that the L\'evy specification provides robust improvements over the Gaussian OU benchmark in both parameter identification and extreme-event tracking, and that these improvements are largely insensitive to the specific choice of $\alpha$. For our baseline results, we adopt $\alpha = 1.2$, as it delivers the strongest alignment with PWT benchmarks and the most accurate crisis tracking for Argentina. The optimal choice of $\alpha$, however, is left for future research.
\subsection{Cross-country evidence}
		
\subsubsection{Motivation and country selection}
		
		The preceding sections establish that the L\'evy-driven Solow model delivers coherent parameter estimates and superior crisis-tracking for Argentina. To assess whether these gains generalize beyond Argentina's turbulent environment, we apply the framework to two additional economies: \textbf{Colombia}, a middle-income emerging economy with moderate volatility, and the \textbf{United States}, a low-volatility advanced economy. This design tests whether the L\'evy specification offers advantages where extreme shocks are rare.
		
\subsubsection{Parameter estimates across countries}
		
		\begin{table}[htbp]
			\centering
			\caption{Parameter Estimates Across Countries: L\'evy vs. Gaussian OU}
			\begin{tabular}{lcccccccccc}
				\toprule
				Country & Method & $\hat{\alpha}_k$ & $\hat{s}_1$ & $\hat{s}_2$ & $\hat{\eta}$ & $\hat{\mu}$ & $\hat{\nu}$ & $\hat{\gamma}$ & $\hat{a}_1$ & $\hat{\sigma}_1$ \\
				\midrule
				Argentina & L\'evy (our method) & --- & 0.076 & 0.155 & \textbf{0.050} & 0.980 & 0.879 & 25.81 & 0.145 & 0.302 \\
				(2008--2023) & OU (Brannan 2019) & 0.378 & 0.135 & 0.429 & 0.360 & 0.082 & 0.989 & 20.59 & 0.062 & 0.321 \\
				& PWT reference & 0.415-0.661
				& 0.110 & 0.180 & --- & & & & & \\
				\midrule
				Colombia & L\'evy (our method) & --- & 0.160 & 0.192 & \textbf{0.052} & 0.123 & 0.261 & 9.86 & 0.129 & 0.303 \\
				(2005--2023) & OU (Brannan 2019) & 0.436 & 0.317 & 0.329 & 0.482 & 0.791 & 0.346 & 19.12 & 0.232 & 0.057 \\
				& PWT reference & 0.462-0.555 & 0.133 & 0.231 & --- & & & & & \\
				\midrule
				United States & L\'evy (our method) & --- & 0.183 & 0.254 & \textbf{0.050} & 0.494 & 0.217 & 6.03 & 0.099 & 0.318 \\
				(2005--2023) & OU (Brannan 2019) & 0.279 & 0.299 & 0.325 & 0.460 & 0.812 & 0.302 & 19.86 & 0.215 & 0.031 \\
				& PWT reference & 0.394-0.432 & 0.200 & 0.276 & --- & & & & & \\
				\bottomrule
			\end{tabular}
			\label{tab:cross_country}
		\end{table}
		
		Table~\ref{tab:cross_country} reports the estimated parameters for all three countries under both the L\'evy specification and the Gaussian OU benchmark, following the format of Tables~1--2. Three findings merit emphasis.
		
		First, the L\'evy specification delivers parameter estimates that are systematically closer to the PWT external benchmarks than those obtained from the Gaussian OU model. This advantage is most pronounced for Argentina---the most volatile economy in our sample---where the L\'evy estimates of \(s_1\), \(s_2\), and \(\eta\) align with the PWT references substantially better than their OU counterparts. For Colombia, a moderately volatile emerging economy, the improvement remains clear though more moderate. For the United States, where business cycles are smooth and extreme shocks are rare, the two specifications yield comparable estimates. This graded pattern is precisely what the theory predicts: the Gaussian OU model, lacking the ability to generate discontinuous jumps, systematically misattributes the persistence of extreme shocks to structural parameters, thereby distorting \(s_1\), \(s_2\), and \(\eta\) most severely in environments where such shocks occur frequently. Our L\'evy framework, by explicitly modeling jumps, isolates the genuine structural components from transient fluctuations, yielding estimates that remain anchored to their fundamental values regardless of the economy's volatility level.
		
		Second, as a corollary to this finding, the capital adjustment speed \(\hat{\eta}\) estimated under the L\'evy specification exhibits remarkable stability across countries---approximately 0.050 for Argentina, 0.052 for Colombia, and 0.050 for the United States. This invariance is consistent with the interpretation that \(\eta\) captures the genuine physical speed of capital adjustment---governed by construction cycles, equipment delivery, installation, and capacity ramp-up, rather than a composite parameter contaminated by shock persistence. Under the Gaussian OU benchmark, by contrast, \(\hat{\eta}\) is inflated by a factor of 7 to 10 (0.360--0.482), confirming that the diffusion specification systematically distorts structural identification when extreme shocks are present. While this stability supports the coherence of our framework, it is ultimately a consequence of the broader identification gains documented above.
		
		Third, the procyclical pattern (\(s_1 < s_2\)) persists across all countries, with magnitudes that vary systematically with economic volatility. The innovation scale \(\hat{\sigma}_1\) also reflects the same underlying logic: it remains stable around 0.30--0.32 under the L\'evy specification, while the Gaussian OU benchmark produces distorted estimates that vary erratically across countries. These patterns are consistent with the interpretation that the L\'evy framework separates transient fluctuations from structural dynamics more cleanly than its Gaussian counterpart.
		
		Taken together, these results demonstrate that the L\'evy framework delivers parameter estimates that systematically outperform the Gaussian OU benchmark against external PWT references. The advantage is largest for the most volatile economies, confirming that the methodological innovations introduced in this work are most valuable precisely where they are most needed.
		
\subsubsection{Predictive performance}
		
\begin{figure}[htbp]
			\centering
			\includegraphics[width=0.9\textwidth]{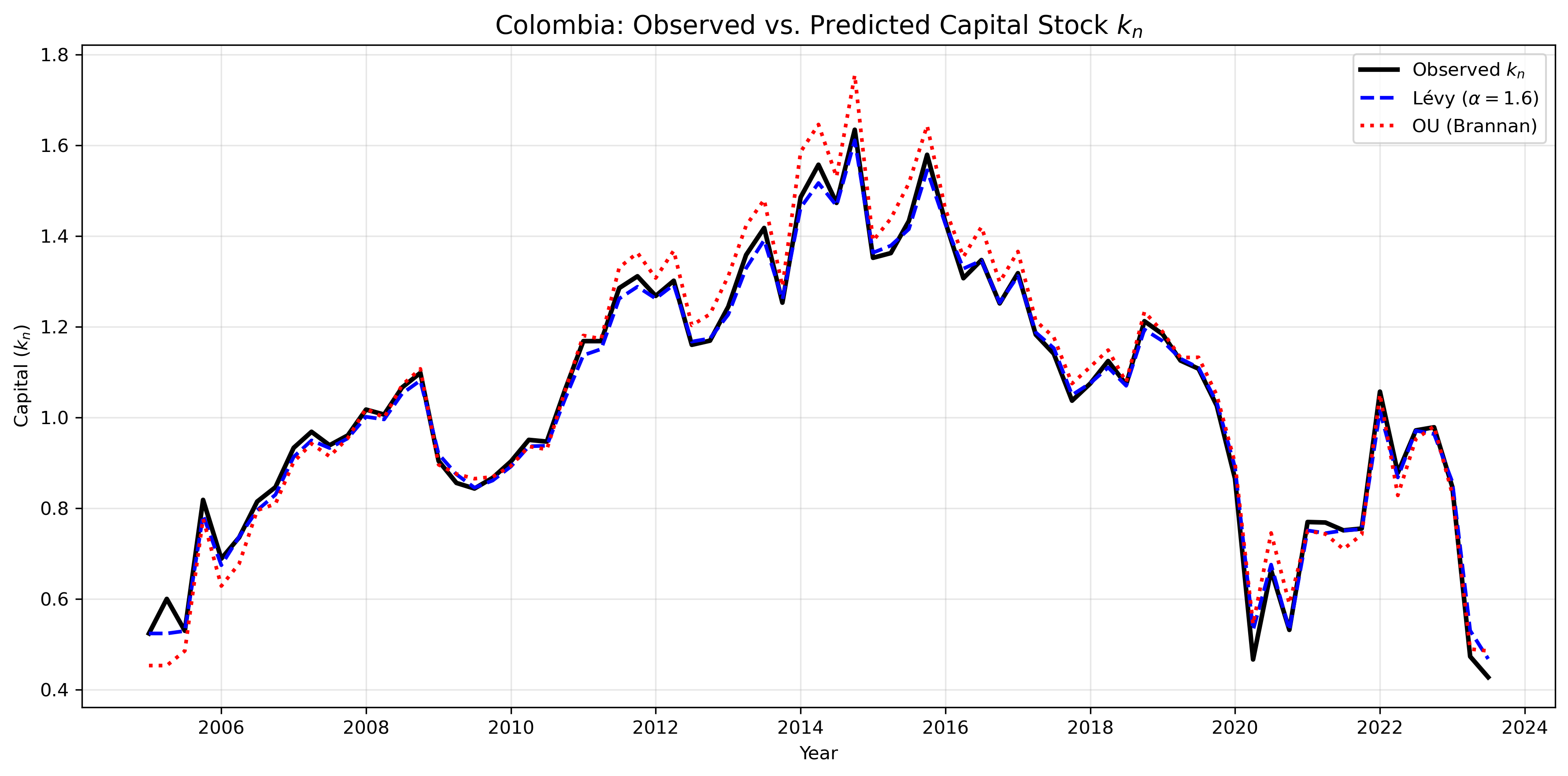}
			\caption{Colombia: Observed vs. Predicted Capital Stock $k_n$}
			\label{fig:colombia_pred}
		\end{figure}

\begin{figure}[htbp]
			\centering
			\includegraphics[width=0.9\textwidth]{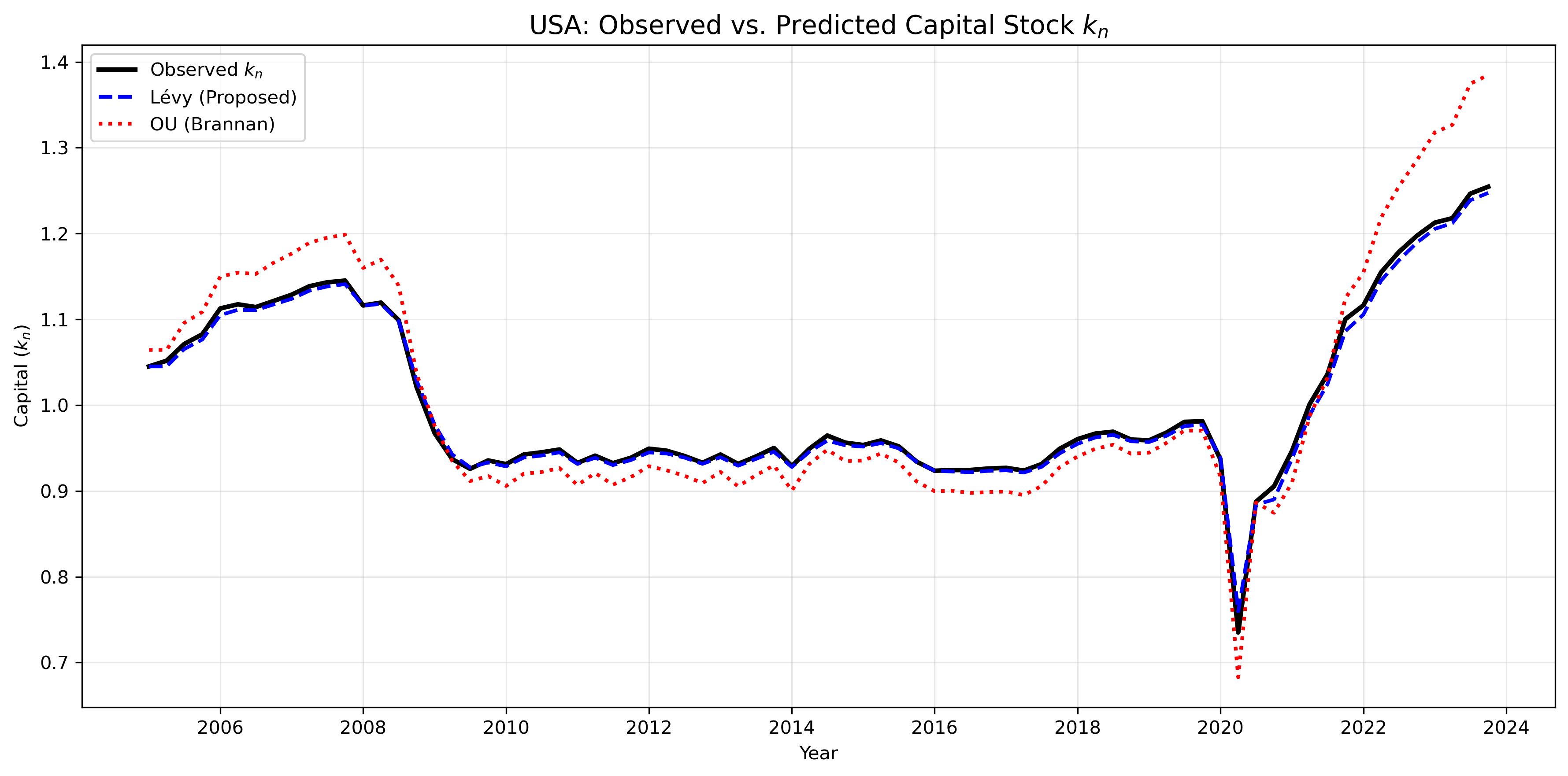}
			\caption{United States: Observed vs. Predicted Capital Stock $k_n$}
			\label{fig:usa_pred}
		\end{figure}

		Figures \ref{fig:colombia_pred} and \ref{fig:usa_pred} compare one-step-ahead predictions from the L\'evy and Gaussian OU models for Colombia and the United States, respectively.
		
		For Colombia (Figure \ref{fig:colombia_pred}), the L\'evy prediction (dashed blue line) tracks the observed series (black solid line) closely across the full sample. During the COVID-19 recession of 2020---the most extreme shock in Colombia's recent history---the L\'evy model captures the sharp trough with substantially greater accuracy than the Gaussian OU benchmark (red dashed), which smooths over the downturn. This mirrors the Argentine evidence from Figures \ref{fig:k_pred_1}--\ref{fig:k_pred_2}.
		
		For the United States (Figure \ref{fig:usa_pred}), the two predictions are broadly similar over most of the sample period. However, during the 2007--2008 slowdown preceding the global financial crisis and again in 2022--2023, the L\'evy prediction tracks the observed series more closely than the Gaussian OU benchmark, which exhibits smoother but less responsive dynamics. The relative proximity of the two predictions for the remainder of the sample reflects the generally mild nature of U.S. business cycles, where Gaussian approximations are adequate for most point-forecasting purposes. Nevertheless, the L\'evy specification performs equivalently or better during tranquil periods while delivering structurally coherent parameters---a combination that the Gaussian benchmark cannot achieve. This result confirms that the L\'evy framework is a genuine generalization of the Gaussian model: it converges to the Gaussian solution when extreme shocks are absent, yet retains the flexibility to capture even moderate fluctuations when they arise.

\section{Conclusion and outlook}
		
		This work has developed and empirically implemented a Solow-type growth model driven by symmetric $\alpha$-stable L\'evy shocks with time-varying capital elasticity. We synthesize the main contributions, discuss their implications, and identify avenues for future research.
		
		\subsection{Main contributions}

This paper has established a Solow-type framework that accommodates $\alpha$-stable L\'evy shocks and time-varying capital elasticity. Theoretically, we provide a complete probabilistic characterization of the capital deviation process---stationary distribution, conditional characteristic function, and an integral representation that clearly separates investment gestation lags from endogenous feedback. Methodologically, we devise an estimation strategy that relies only on finite moments (mean absolute error and mean matching), bypassing the non-existence of variance in L\'evy-driven data. Empirically, our framework yields structurally invariant parameters ($\eta \approx 0.05$) across Argentina, Colombia, and the United States, and substantially improves crisis-period tracking for volatile economies without sacrificing performance in tranquil environments. These results demonstrate that the L\'evy specification is not an ad hoc deviation from the Gaussian benchmark; it converges to the Gaussian solution when extreme shocks are absent and captures them more accurately than the Gaussian benchmark when they occur. This dual property confirms the framework's theoretical coherence and empirical robustness.

		\subsection{ Policy implications}
		
		Four implications emerge from our analysis.
		
		First, policymakers in high-volatility economies should recognize that extreme shocks are recurrent, not exceptional. Conventional Gaussian-based tools systematically underestimate crisis probabilities and magnitudes. The L\'evy framework developed in this paper provides a practical toolkit for more reliable growth projections, fiscal planning, and external vulnerability assessments under these conditions.
		
		Second, the capital elasticity of output, \(\alpha_k\), is not a purely technological parameter. It reflects the distribution of national income between capital and labor, shaped by a country's political economy, labor market institutions, and social welfare systems. The substantial year-to-year variation in \(\alpha_k^t\) for Argentina signals instability in its distributional arrangements. Stabilizing \(\alpha_k\) ultimately requires a stable institutional framework for income distribution, not merely factor accumulation policies.
		
		Third, the cross-country stability of \(\eta \approx 0.05\) implies that capital adjustment is inherently slow, with a half-life of approximately 3.5 years. Demand-management policies can at best temporarily offset deviations from the balanced growth path; accelerating long-run capital accumulation requires structural reforms that directly reduce adjustment frictions.
		
		Fourth, the saving behavior recovered under the L\'evy specification reveals that precautionary motives intensify during crises in volatile economies, dampening the consumption response to fiscal stimulus. Countercyclical policy should therefore be complemented with measures that directly address uncertainty---credible policy commitments, exchange rate stabilization, and social safety nets.
		
		\subsection{Future challenges}
		
		This section discusses the limitations of our study and identifies directions for future work.
		
		\paragraph{Data and measurement} We treat the capital elasticity \(\alpha_k^t\) as a known input calibrated from external labor-share data. This approach reduces the number of parameters to be estimated and improves identification, but it also introduces measurement error---particularly in economies with large informal sectors, where official labor-share statistics may be systematically biased. This limitation, however, also points to a deeper policy insight: \(\alpha_k\) is not a purely technological parameter but a reflection of income distribution. If the calibration of \(\alpha_k^t\) could be elevated from a passive external input to a dynamically adjustable policy variable, the model could offer quantitative guidance for structural reforms. Future work could estimate \(\alpha_k^t\) jointly with other parameters or cross-validate macro-calibrated results with micro-level data to better identify the institutional drivers of variation in \(\alpha_k^t\).

			\paragraph{Model specification} We assume symmetric \(\alpha\)-stable shocks, implying identical tail behavior for positive and negative innovations. Evidence from emerging economies suggests that negative shocks---sudden stops, currency crises---are typically more severe and frequent than positive ones. Extending the framework to skewed L\'evy specifications would capture this asymmetry and improve the model's empirical fit for such economies. Relatedly, while our time-varying \(\alpha_k^t\) strategy provides a more flexible production structure, it remains a passive calibration of each country's actual distributional regime. Allowing \(\alpha_k^t\) to vary endogenously across countries---or to evolve according to country-specific institutional characteristics---could further enhance the framework's predictive power and its applicability to cross-country forecasting.
		
			\paragraph{Methodological extensions} The conditional characteristic function derived in Subsection 2.4 provides a complete probabilistic characterization of future capital deviations. While we have focused on point forecasts via the conditional mean, the full distribution could be exploited to construct predictive intervals and density forecasts---tools of direct value for stress testing, macroprudential policy design, and extreme scenario analysis.
		
			\paragraph{External validity and open-economy extensions} Our cross-country evidence covers three economies with distinct volatility regimes, but the sample remains limited. Expanding the analysis to additional emerging and advanced economies would strengthen the generality of our findings. In addition, our analysis is conducted within a closed-economy framework. Extending the model to incorporate international capital flows, exchange rate dynamics, and terms-of-trade shocks would enhance its relevance for small open economies like Argentina and Colombia, where external vulnerabilities are a primary source of macroeconomic volatility.

\bigskip
\noindent\textbf{Data availability}

The numerical algorithms and source code that support the findings of this study are available from the corresponding author upon reasonable request.

\bigskip
\noindent\textbf{Declaration of competing interest}

No author associated with this paper has disclosed any potential or pertinent conflicts which
may be perceived to have impending conflict with this work.

\bigskip
\noindent\textbf{Acknowledgments}

 This work was supported by the Guangdong Basic and Applied Basic Research Foundation (Grant No. 2025A1515012560), the Guangdong Introduction Program (Grant No. 2023QN10X753) and National
Foreign Experts Program (Grant No. 111001819820258003).

\end{document}